\documentclass[12pt, amssymb, amsmath]{article}

\begin{document}

\renewcommand{\thesection}{\arabic{section}}
\renewcommand{\thesubsection}{\thesection.\arabic{subsection}}


\newtheorem{lem}{Lemma}[section]
\newtheorem{propo}[lem]{Proposition}
\newtheorem{theo}[lem]{Theorem}
\newtheorem{rema}[lem]{Remark}
\newtheorem{remas}[lem]{Remarks}
\newtheorem{coro}[lem]{Corollary}
\newtheorem{defin}[lem]{Definition}
\newtheorem{hypo}[lem]{Hypoth\`ese}
\newtheorem{exem}[lem]{Exemple}
\newtheorem{conj}[lem]{Conjecture}

\newcommand{\real}{{\bf R}}
\newcommand{\reall}{{\bf R}^{p}}
\newcommand{\realll}{{\bf R}^{k}}
\newcommand{\mreal}{M\times{\bf R}^{p}}
\newcommand{\mreall}{M\times{\bf R}^{2k}}
\newcommand{\ri}{\rightarrow}
\newcommand{\tast}{T^{\ast}M}
\newcommand{\tastl}{T^{\ast}L}
\newcommand{\tastilde}{T^{\ast}\widetilde{M}}
\newcommand{\cd}{D_{\bullet}}
\newcommand{\cdk}{D_{\bullet-k}}
\newcommand{\cfx}{C_{\bullet}(f,\xi )}
\newcommand{\cald}{{\cal D}(M)}
\newcommand{\cm}{C_{\bullet}(M)}
\newcommand{\cmk}{C_{\bullet-k}(M)}
\newcommand{\cov}{\widetilde{M}}
\newcommand{\calm}{{\cal M}}
\newcommand{\cala}{{\cal A}}
\newcommand{\calc}{{\cal C}}
\newcommand{\call}{{\cal L}}
\newcommand{\caly}{{\cal Y}}
\newcommand{\gol}{\, \, \, \, \, \, \, \, }
\newcommand{\hatf}{\hat{f}}
\newcommand{\hatx}{\hat{x}}
\newcommand{\hatv}{\hat{v}}
\newcommand{\mn}{M^{n}}
\newcommand{\barm}{\bar{M}}
\newcommand{\bfs}{{\bf S}}
\newcommand{\piu}{\pi_{1}}
\newcommand{\pii}{\pi_{i}(M)}
\newcommand{\omi}{\Omega(L_{0},L_{1})}
\newcommand{\omibar}{\Omega(\bar{L}_{0},\bar{L}_{1})}
\newcommand{\omil}{\Omega(L,L)}
\newcommand{\omibarl}{\Omega(\bar{L},\bar{L})}
\newcommand{\la}{^}
\newcommand{\bfc}{{\bf C}}
\newcommand{\covl}{\widetilde{L}}
\newcommand{\barl}{\bar{L}}

\begin{center}
{\large\bf  FLOER HOMOLOGY ON THE UNIVERSAL COVER, A PROOF OF AUDIN'S CONJECTURE
 AND OTHER CONSTRAINTS ON LAGRANGIAN SUBMANIFOLDS}\\
\vspace{.3in}
\noindent Mihai DAMIAN\footnote{Supported by ANR project "Floer Power" ANR-08-BLAN-0291-03}\\
\noindent Universit\'e de  Strasbourg\\
IRMA, 7, rue Ren\'e Descartes,\\
67 084 STRASBOURG\\
e-mail: damian@math.u-strasbg.fr\\

\end{center}
\vspace{.4in}

\noindent{\bf Abstract\,:}\, We establish a new version of Floer homology for monotone Lagrangian
embeddings in symplectic manifolds. As applications, we get assertions 
for (monotone) Lagrangian submanifolds $L\hookrightarrow
M$ which are
displaceable through Hamiltonian isotopies (this happens for instance when $M={\bf C}\la{n}$). We show
that when $L$ is aspherical, or more generally when the homology of its universal cover vanishes in odd
degrees,  its Maslov number $N_{L}$ equals $2$. We also give topological characterisations of Lagrangians
$L\hookrightarrow M$ with maximal Maslov number: when $N_{L}=dim(L)+1$ then $L$ is homeomorphic to a sphere; 
when
$N_{L}=n\geq 6$ then $L$ fibers over the circle and the fiber is homeomorphic to a sphere. A
consequence is that any exact Lagrangian in $T\la{\ast}{\bf S}\la{2k+1}$ whose
Maslov class is zero is homeomorphic to $\bfs\la{2k+1}$. 
\vspace{.2in}

\section{Introduction and main results} 

\subsection{Preliminaries}

  Let  $(M^{2n}, \omega)$ be a symplectic manifold. A submanifold $L^{n}$ of $M$ is called
Lagrangian if the restriction of $\omega$ on $L$ vanishes. Throughout this paper all 
 symplectic
manifolds are assumed to be either  closed or convex at infinity and all 
Lagrangian submanifolds are assumed to be closed and connected.  One of the fundamental 
questions
in symplectic geometry is the following:
\begin{quote} What properties has to satisfy a closed manifold $L$ in order to admit a
Lagrangian embedding into a given symplectic manifold~$M$?
\end{quote} 

 This question is still widely open even in the case of $M=\bfc\la{n}$. 
The results of the present paper  concern symplectic manifolds such as $M=\bfc\la{n}$,
$M={\bf CP}\la{n}$, or $M=T\la{\ast}K$ for $K$ is closed, all of them being endowed with
their standard symplectic form.  
We establish new topological constraints on Lagrangian submanifolds
$L\subset M$ which are {\it monotone} or {\it exact}. These notions are defined using the following 
two
morphisms related
 to a given  Lagrangian $L$. The morphism
$I_{\omega}:\pi_{2}(M,L)\ri\real$ is defined by:
$$I_{\omega}(A)=\int_{A}\omega.$$
In order to define the morphism $I_{\mu} :\pi_{2}(M,L)\ri{\bf Z}$, pick a smooth map of
pairs $w:(D\la{2},\partial D\la{2})\ri(M,L)$ in the class $A\in\pi_{2}(M,L)$. There is an
unique trivialisation (up to homotopy) of the pull-back $w\la{\ast}TM\, \approx\,
D\la{2}\times{\bf C}\la{n}$ as a symplectic vector bundle. This gives a map
$\alpha_{w}$ from ${\bf S}\la{1} =\partial D\la{2}$ to $\Lambda({\bf C}\la{n})$ - the set of
Lagrangian planes in ${\bf C}\la{n}$. On this space there is a well-known Maslov class
$\mu\in H\la{1}(\Lambda({\bf C}\la{n}),{\bf Z})$ (see \cite{Ar2}), so that one can define 
$$I_{\mu}(A)\, =\, \mu(\alpha)\in{\bf Z}.$$

\begin{defin} A Lagrangian submanifold $L\subset M$ is called {\it weakly exact} if the
morphism $I_{\omega}$ vanishes. It is called {\it exact} if $\omega=d\lambda$ and the
restriction $\lambda|_{L}$ is an exact one-form. 

A Lagrangian submanifold is called {\it monotone} if there is a constant $\tau>0$ such 
that 
$$I_{\omega}\, =\, \tau I_{\mu}.$$

\end{defin}

By definition, only exact symplectic
manifolds admit exact Lagrangian submanifolds. It is not obvious, but still true, 
that  monotone Lagrangian
 submanifolds only exist in monotone symplectic manifolds (i.e. in symplectic manifolds 
in which the
morphism defined on $\pi_{2}(M)$ by the symplectic form  is a positive multiple of 
the morphism defined
by the first Chern class).
   Many authors studied monotone and exact Lagrangians and found various obstructions 
to the existence
of such embeddings. A celebrated result of M. Gromov asserts:

\begin{theo}\label{Gromov} \cite{G} There is no weakly exact Lagrangian embedding 
$L\subset
{\bf C}\la{n}$.  
\end{theo}

The results on the obstructions to the existence of
 monotone Lagrangian submanifolds  mostly concern their Maslov number. This
number, denoted by $N_{L}$, is defined as the minimal positive integer which is in the
image of $I_{\mu}$. In  1996 Y.-G. Oh established the following inequality \cite{Oh}, 
improving thus
previous results of C. Viterbo \cite{Vi} and L. Polterovich \cite{Polt1}, \cite{Polt2}:

\begin{theo}\label{oh} For any monotone Lagrangian submanifold $L\subset{\bf C}\la{n}$ 
we have~:
$$1\, \leq\, N_{L}\, \leq n.$$
\end{theo}

  These bounds turn out to be sharp. Indeed, L. Polterovich gave in \cite{Polt2}
 an example of a monotone
Lagrangian $L\subset{\bf C}\la{n}$ which satisfies $N_{L}=n$. 

Note that both Gromov's and Oh's result can be stated in the  more general case of symplectic
manifolds $M$ which are convex at infinity and have the property that any compact subset 
is displaceable through a
Hamiltonian isotopy. This means that for any compact $K\subset M$ there is a Hamiltonian
isotopy $(\phi_{t})_{t\in[0,1]}$ such that $\phi_{1}(K)\cap K = \emptyset$. Symplectic 
manifolds of the
form ${\bf C}\times W$,
 or subcritical Stein manifolds \cite{BiC} satisfy  this assumption. 

However, in this more general case the statement of Oh's result is slightly different:

\begin{theo}\label{oh2} Let $M$ be a symplectic manifold in which every compact 
subset is displaceable
through a Hamiltonian isotopy. For any monotone Lagrangian submanifold $L\subset M$ we
 have:
$$1\, \leq\, N_{L}\, \leq n+1$$
and if $N_{L}=n+1$, then $L$ is a ${\bf Z/2}$-homology sphere.
\end{theo}

  Actually,  more recent results of K. Fukaya, Y.-G. Oh, H. Ohta and K. Ono \cite{FO3}, 
imply that for
$N_{L}=n+1$ the Lagrangian $L$ is a $\bf Z$-homology sphere in the statement above. 

\subsection{Main results} 

Our results about monotone Lagrangian submanifolds are of two types. First we show that 
under some
topological assumptions on $L$ we have $N_{L}=2$. Then, we study the topology of  monotone 
Lagrangian
submanifolds with maximal Maslov number $N_{L}=n+1$ or $N_{L}=n$.
 Here are the statements:

\begin{theo}\label{main1} Let $M$ be a monotone symplectic manifold which has the 
property that 
 any compact subset
is displaceable through a Hamiltonian isotopy. Let  $L\subset M$ be a monotone Lagrangian
submanifold.\\
(a) If $L$ is aspherical (i.e. $L$ is an Eilenberg Mc Lane space $K(\pi_{1}(L),1)$), 
then $N_{L}=2$ if $L$ is
orientable and $N_{L}\in\{1, 2\}$ if $L$ is not orientable.\\
(b) Denote by $\covl$ the universal cover of $L$. If $L$ is orientable and has the property 
$$H_{2i+1}(\covl,{\bf Z/2})\, =\, 0$$ for any integer $i$ then $N_{L}=2$. \\
(c) Moreover, for any
almost complex structure $J$ which is compatible with the symplectic form, a Lagrangian 
$L$ which satisfies  the condition (b)  has the 
property
that through every $p\in L$ there is  a $J$-holomorphic disk $w:(D,\partial D)\ri (M,L)$ such that:
\begin{itemize}
\item   The Maslov
index $\mu(w)$ equals $2$.
\item  $p\in w(\partial D)$.
\item $w(\partial D)$ is non zero in $\piu(L)$.
\end{itemize} 
\end{theo}

\vspace{.2in}

\noindent{\bf Remarks} \\
\noindent 1. Part (a) of Theorem \ref{main1} was proved by K. Fukaya 
 for general aspherical Lagrangian
submanifolds, but under the additional hypothesis that $L$ is orientable and 
 relatively spin (\cite{Fu1}, Th. 12.2). In
the case where $L$ is a torus the statement was conjectured by M. Audin \cite{Au}. 
In this particular
case, many  results were previously obtained by L. Polterovich, C. Viterbo, Y.-G. Oh, 
Y. Eliashberg, 
P. Biran, K. Cieliebak
and L. Buhovsky. \\
2. Part (b) of the statement above was proved by K. Fukaya in the case $L={\bf S}
\la{1}\times{\bf
S}\la{2m}$ without any monotonicity assumption (\cite{Fu1}, Th. 13.1).
 However, our result applies to many more general 
examples, such as arbitrary products of tori (or other orientable aspherical manifolds) and 
complex projective spaces, 
even-dimensional
spheres, etc. \\
3. Many results related to  part (c) of the theorem can be found in the paper \cite{BiCo3} of P. Biran
and O. Cornea. Using their terminology, the Lagrangian $L$ should be called
 {\it uniruled of type $(0,1)$ and order
$2$.} \\

Using the ideas of P. Biran \cite{Bi}, we obtain the following corollary on the monotone 
Lagrangian
submanifolds in the complex projective space:

\begin{theo}\label{coro1} Let $W$ be a symplectic manifold such that 
$M={\bf CP}\la{n}\times W$
 is
monotone (for instance this holds for $\pi_{2}(W)=0$ or for $W={\bf CP}\la{n}$). 
Let $L\subset M$
 be
a monotone Lagrangian submanifold which is
aspherical. Then $N_{L}=2$ if $L$ is orientable and $N_{L}\in\{1,2\}$ if $L$ is not orientable.\\ 
\end{theo}

For  spin Lagrangian submanifolds and $W=point$ this result
 was also proved by K. Fukaya 
in \cite{Fu1} without any monotonicity
assumption. The result is still true in the more general situation where ${\bf CP}\la{n}$ 
is replaced
by a symplectic manifold which arises  as a hypersurface in a subcritical polarisation. 
These manifolds
were studied in \cite{BiC} by P. Biran and K. Cieliebak. \\ 

Our next result is a topological characterisation of monotone Lagrangian submanifolds with 
maximal
Maslov number.

\begin{theo}\label{main2} Let $M$ be a monotone symplectic manifold of dimension $2n$,  
which has the property that 
 any compact subset
is displaceable through a Hamiltonian isotopy. Let  $L\subset M$ be a monotone Lagrangian
submanifold.\\
(a) Suppose that $N_{L}=n+1$ and $n\geq 2$. Then $n$ is odd and $L$ is homeomorphic to a 
$n$-sphere. \\
(b) Suppose that $N_{L}=n$ and $n\geq 3$.

If $n$ is odd then $\pi_{1}(L)$ 
has an infinite cyclic group $G\approx {\bf Z}$ of finite index. If moreover $M$ is an exact 
symplectic
manifold then  there is an exact sequence of groups  $$0\ri K\ri \pi_{1}(L)\, \ri\, {\bf Z}\ri 0,$$
where $K$ is finite and has odd order.

 If $n$ is even then  $\pi_{1}(L)\, \approx\, {\bf Z}$. If  moreover $n\, \geq 6$,
then there is a fibration of $L$ over the circle ${\bf S}\la{1}$ whose fiber is 
homeomorphic 
to the
$(n-1)$-sphere.
\end{theo}

\vspace{.2in}

\noindent{\bf Remarks}\\
 1. There are examples of monotone Lagrangian submanifolds satisfying the
hypothesis on the Maslov number above. Indeed the embedding of ${\bf
S}\la{2k+1}$ into ${\bf CP}\la{k}\times \bfc\la{k+1}$ given by 
$$z\, \mapsto\, ([z],\bar{z})$$ is monotone, Lagrangian and its Maslov number is $2k+2$.
 This example
is due to M. Audin, F. Lalonde and L. Polterovich \cite{ALP}. An example of a monotone 
Lagrangian 
embedding
 ${\bfs}\la{1}\times\bfs\la{2k-1}\, \subset \bfc\la{2k}$ whose Maslov number equals $2k$
 was
constructed by L. Polterovich in \cite{Polt2}.\\
2. In \cite{Gad} A. Gadbled established topological constrains on
 monotone Lagrangian submanifolds in cotangent bundles which have a large Maslov number (which implies
that they are not displaceable throgh Hamiltonian isotopies).\\
 
We prove the following corollaries of this theorem:

\begin{theo}\label{projective} (a) Let $X$ be a symplectic manifold of dimension $2n+2$ with
$\pi_{2}(X)=0$. Let $L\subset {\bf CP}\la{n}\times X$ be a Lagrangian 
submanifold such that $H_{1}(L,
{\bf Z})=0$. Then $L$ is homeomorphic to $\bfs\la{2n+1}$. \\
(b) Let $L\subset {\bf CP}\la{n}\times{\bf CP}\la{n}$ be a Lagrangian 
submanifold such that $H_{1}(L,
{\bf Z})$ vanishes. Then $L$ is simply connected and there is a circle 
fibration $\bfs\la{2n+1}\ri L$. \\
(c)
Let $L\subset {\bf CP}\la{n}$ be a Lagrangian submanifold 
with $2H\la{1}(L, {\bf Z})\, =\, 0$ . Then  if $n$ is odd we have $\piu(L)={\bf Z/2}$ and
the universal cover of $L$ is homeomorphic to $\bfs\la{n}$. 
\end{theo} 

 In  \cite{BiCo} O. Cornea and P. Biran asked whether a Lagrangian as in Theorem \ref{projective}.(c)
is diffeomorphic (or homeomorphic) to ${\bf RP}\la{n}$. Our result goes in this direction but
we do not know whether its conclusion implies that $L$ is homeomorphic to the projective
space.  
In the mentioned
paper Biran and Cornea proved that under the hypothesis of Theorem
\ref{projective} the cohomology ring  (with ${\bf
Z/2}$-coefficients)
 of $L$
is isomorphic to the cohomology  ring of ${\bf RP}\la{n}$.  Similar results 
 were previously obtained by
P. Biran \cite{Bi} and P. Seidel \cite{Sei}.  Statement (a) 
generalizes Th. B of \cite{Bi} (asserting that $L$ is a
homology sphere). Statement (b) generalizes Th. C of \cite{Bi} 
(which asserts that $L$ has the homology of ${\bf
CP}\la{n}$).  

\begin{theo}\label{sphere} (a)  Let $L\subset T\la{\ast}\bfs\la{2k+1}$ be an exact 
Lagrangian 
submanifold
with vanishing Maslov class. Then $L$ is homeomorphic to $\bfs\la{2k+1}$. \\
(b)  Let $K\la{2k+1}$  be a manifold whose universal cover is  $\bfs\la{2k+1}$.
Let $L\subset T\la{\ast}K$ be an exact Lagragian submanifold
with vanishing Maslov class. Then the  universal cover 
$\covl$ is homeomorphic to $\bfs\la{2k+1}$ (in particular
$\piu(L)$ is finite). For instance, when 
 $K={\bf RP}\la{2k+1}$ then $\piu(L)={\bf Z/2Z}$ and 
 $L$ is double covered by (a manifold homeomorphic to) ${\bfs}\la{2k+1}$. \\
\end{theo}

This result gives an answer  in the case $K=\bfs\la{2k+1}$ (under the
hypothesis of the vanishing Maslov class) to an open question
raised by V.I. Arnold  \cite{Ar}: \begin{quote} Is an exact submanifold 
$L\subset T\la{\ast}K$ homeomorphic to
$K$ ? 
\end{quote}

Actually, Arnold asks whether $L$ is Hamiltonian isotopic 
to the zero section, but this latter
question seems out of reach, except for the case $dim(N)=2$ 
(see \cite{Hi}  for related
results). In the general case the most striking result was obtained by 
K. Fukaya, P. Seidel and I.
Smith, who proved in \cite{FSS1}, \cite{FSS2} that, when $L$ is
relatively spin  with vanishing
Maslov class, its cohomology is isomorphic to the cohomology of $K$. 
For the case $K=\bfs\la{m}$
similar results were previously obtained by P. Seidel \cite{Sei} and L. Buhovsky \cite{Buho}.

\subsection{Idea of the proof}

 Recall first the definition of Floer homology.
 Let $L\subset M$ be a Lagrangian submanifold which is monotone with $N_{L}\geq 2$,  
or weakly
 exact. Consider a Hamiltonian isotopy
$(\phi_{t})$ defined by a time-dependent Hamiltonian $H:[0,1]\times M\ri\real$
 and a time-dependent almost complex structure $J$  which is compatible to the symplectic form $\omega$.
 For a generic
choice of the couple $(H,J)$, A.
Floer associated to these data a complex $(C_{\bullet}(H), \partial_{J})$ whose homology  
does not depend on $(H,J)$ \cite{F1}, \cite{F2}, \cite{F3}. The Floer complex is 
free over 
${\bf Z/2}$, 
spanned by the intersections $L\cap\phi_{1}(L)$, which are
supposed to be transverse. Its differential $\partial_{J}$ is defined by counting 
the {\it isolated} 
holomorphic
strips $$v:\real\times[0,1]\ri M,$$
with boundary in $L\cup \phi_{1}(L)$ (more precisely $v(\real\times\{i\}) \subset 
\phi_{i}(L)$
 for
$i=0, \, 1$) and joining intersection points $x,y \, \in\, L\cap\phi_{1}(L)$, which 
means that 
$$\lim_{s\ri -\infty}v(s,t) = x \, \, \, \mbox{and} \, \, \, 
\lim_{s\ri +\infty}v(s,t) = y.$$
If $n(x,y)$ is the number modulo $2$ of such curves then 
$$\partial_{J}(x)\, =\, \sum_{y\in L\cap\phi_{1}(L)}n(x,y)y.$$
Note that if $L$ is orientable and 
relatively spin, (meaning that  the second 
Stiefel-Whitney class $w_{2}(L)$
lies in the image of $H\la{2}(M,{\bf Z/2})\ri H\la{2}(L, {\bf Z/2})$), then the 
whole theory
 works for
integer coefficients \cite{FO3}. 

A relation between the Floer homology $HF(L)$ and the usual homology can be established. 
At this end, one
should remark that given a Morse function $f:L\ri \real$ which is sufficiently 
$\calc\la{2}$-small, 
its
graph $$\{(df_{q}, q)\, |\, q\in L\}\, \, \subset\, \, T\la{\ast}L$$
can be embedded in $M$ via a Weinstein neighborhood $U(L)\subset M$ and the 
intersection points 
$L\cap L_{f}$ correspond to
the critical points of $f$. Moreover, for a good choice of the almost complex 
structure $J$, the 
application 
$$v(s,t) \, \mapsto\, v(s,0)$$ defines a one-to-one correspondence between the 
holomorphic strips 
joining
two intersection points $x, y$ {\it which lie in $U(L)$}
 and the flow lines of a vector field on $M$ which is the gradient of $f$ with
respect to some Riemannian metric. So the Morse complex becomes a sub-complex of 
the Floer complex 
in this
case. If $L$ is weakly exact, and $f$ is chosen sufficiently small, one can prove that 
no holomorphic strip leaves
$U(L)$, so that the two complexes are actually isomorphic and thus Floer homology is 
isomorphic to 
usual
homology. In the case where $L$ is monotone, Y.-G. Oh shows in \cite{Oh} that in the 
case of the
particular Hamiltonian isotopy defined by the graph of a small function $f$, the 
Floer differential
decomposes into a sum 
$$\partial_{J} \, =\, \partial_{0}+\partial_{1}+\partial_{2}+\cdots,$$
where $\partial_{0}:C_{k}(f)\ri C_{k-1}(f)$ is the Morse differential and 
$\partial_{l}:C_{k}(f)\ri
C_{k-1+lN_{L}}(f)$ for any integer $l$. By comparing the degrees of the $\partial_{i}$'s 
in the relation
$\partial_{J}\la{2}=0$ one easily sees that $\partial_{1}$ defines an application of 
degree $-1+N_{L}$
 on
the usual homology groups of $L$ and moreover that this application is actually a 
differential. On the
resulting homology groups 
$\partial_{2}$ defines an application of degree $-1+2N_{L}$ which again turns out 
to satisfy
$\partial_{2}\la{2}=0$, and so on... This feature of the Floer differential can be 
formalized in the
 existence
of a spectral sequence which converges to the Floer homology $HF(L)$ and whose first 
page is built 
using the
usual homology groups of $L$ \cite{Oh}, \cite{Bi}. \\

The main idea of our paper is the following. Fix a covering space $\barl\ri L$.
 Given a Morse function on $L$ and
an associated generic gradient, one can build a  free complex, possibly infinite 
dimensional, 
 by lifting its flow 
lines to
$\barl$. The homology of this complex is the usual (singular) homology $H_{\ast}(\barl)$. 

Now consider the Floer complex associated to $L$ and to some generic pair $(H,J)$. 
Look at the 
collection of
paths $s\mapsto v(s,0)\, \subset\, L$ defined by the   holomorphic strips
$v:\real\times[0,1]\ri M$ which are counted in  the Floer differential $\partial_{J}$.  
This collection of  paths obviously contains all the information about the Floer complex; 
one can reconstruct  
it
 by counting  
the paths who 
have the same
endpoints $x$ and $y$ and defining thus a differential on the free module spanned by all these endpoints.

 On the other hand one can lift these paths to the covering space $\barl$. The
question is:
\begin{quote} 
Do the lifted paths define a  complex ? 
\end{quote}

This turns out to be true when $L$ is weakly exact or monotone with $N_{L} \geq 3$. 
But it fails for
$N_{L}=2$, as one can easily see by looking at the example of an embedded circle $L=\bfs\la{1}$ 
in $\bfc$  deformed by a small translation $\phi$, so that $L$ and $\phi(L)$ have
 two
intersection points: in this case the  lifts of the paths defined on $L$ by the  holomorphic strips
of the usual Floer complex 
 do not define a differential  on
 the universal cover of the circle.

When the complex on $\barl$ is defined, we call it {\it lifted Floer complex}.
In order to compute its homology $FH\la{\barl}(L)$ 
we consider a  Hamiltonian isotopy defined by a graph $L_{f}$ of a small
function on $L$. In this case  the differential $\partial_{0}$ clearly  defines the  lifted Morse complex 
described 
above. Therefore, when
it is defined, the lifted Floer homology $FH\la{\barl}(L)$ has analogue features, namely it 
coincides with 
the 
(Morse) homology of  $\barl$ if $L$ is weakly exact and 
it is the limit of  a spectral sequence 
starting from 
$H_{\ast}(\barl)$ when $L$ is monotone. This allows us to prove the claimed results, using properties of the
homology of the covering $\barl$; for instance the fact that in the case of the universal cover
 it vanishes in every nonzero 
degree provided that  $L$ is aspherical. 

This construction will be formalized in the next section. 

\section{The lifted Floer complex}

Let $L\subset M$ be a Lagrangian submanifold. Let $p:\barl\ri L$ be a covering of $L$. 
The elements of a fiber of $p$ are indexed by a
possibly infinite set $I$.
 Let $(\phi_{t})_{t\in [0,1]}$ be a Hamiltonian isotopy of $M$ such that $L$ and 
$\phi_{1}(L)$ are
transverse. For any $x\in L\cap\phi_{1}(L)$ denote by $(x_{i})_{i\in I}$ the elements of
$p\la{-1}(x)$. We prove the following theorem, which is the main ingredient in the 
proof of the
results that we claimed in the preceding section:

\begin{theo}\label{lift} If $L$ is exact or monotone with $N_{L}\geq 3$,  
there exists a free ${\bf Z/2}$-complex $C_{\bullet}$ spanned by
$\bigcup_{x\in L\cap\phi_{1}(L)}\{x_{i}\, |\, i\in I\}$ such that:\\

\begin{itemize}
\item If $L$ is exact then $$H_{\ast}(C_{\bullet})\approx H_{\ast}(\barl,{\bf Z/2}).$$

\item If $L$ is monotone with $N_{L}\geq 3$ then there exist applications $\delta_{1},
\delta_{2}, \ldots, \delta_{k}, \ldots$ with the following properties: 

 (i) $\delta_{1}: H_{\ast}(\barl,{\bf Z/2})\ri  H_{\ast-1+N_{L}}(\barl,{\bf Z/2})$ and
$\delta_{1}\circ\delta_{1}=0$.

 (ii) $\delta_{l}: H_{\ast}(\barl,{\bf Z/2})\ri  H_{\ast-1+lN_{L}}(\barl,{\bf Z/2})$, $l\geq 2$ is
well-defined if $\delta_{m}=0$ for $m=1, \ldots, l-1$ and
$\delta_{l}\circ\delta_{l}=0$.
 
 (iii) If $\delta_{l}=0$ for any $l\geq 1$ then 
$$H_{\ast}(C_{\bullet})\approx H_{\ast}(\barl,{\bf Z/2}).$$

\end{itemize}

When $L$ is orientable and relatively spin, one can replace ${\bf Z}/2$-coefficients with 
${\bf Z}$-coefficients in the statement of this theorem and get in the monotone case
 applications $\delta_{i}$ satisying (i), (ii) above and moreover \\

(iii) If $H_{\ast}(C_{\bullet})$ vanishes there is some $\delta_{l}$ which is not zero and if
$$(l+1)N_{L}>  n+1,$$ then $(H_{\ast}(\barl,{\bf Z}),\delta_{l})$ is acyclic.
\end{theo}

 Before giving the proof of the theorem we briefly remind  the construction of the of the original
 Floer complex $FC_{\bullet}(L, (\phi_{1}(L)))$ over ${\bf Z/2}$. 

\subsection{The usual  Floer complex} 

The
 Floer complex $FC_{\bullet}$ is free, 
spanned by the intersection points $L\cap\phi_{1}(L)$. Denote $L_{t}=\phi_{t}(L)$. 
In order to define the differential of
$FC_{\bullet}$, one has to choose a family of almost complex structures $J=(J_{t})_{t\in[0,1]}$ on $M$ which 
are compatible with
the symplectic form $\omega$ and to define the space of holomorphic strips with bounded energy
$\calm(L_{0}, L_{1})$,  as 
follows: 
$$ \calm(L_{0}, L_{1})\, =\, \left\{v\in\calc^{\infty}(\real\times [0,1],M)
\, \left|\, \begin{array}{c}\frac{\partial v}{\partial
s}+J\frac{\partial v}{\partial t}=0\\ 
\\v(s,0)\in L_{0}, \, v(s,1)\in L_{1}\\
\\
E(v)<+\infty\end{array}\right|\,
\right\}.$$
Here, the energy $E(v)$ is defined by the formula: 
$$ E(v)\, =\, \int_{\real\times[0,1]}\left|\left|\frac{\partial v}{\partial
s}\right|\right|^{2}\, dsdt,$$
the time-dependent norm being defined by the Riemanian metric $\omega(\cdot, J_{t}\cdot)$.
 Then consider  for
$x,y\in L_{0}\cap L_{1}$ the space  
$$\calm(x,y) = \left\{v\in\calm(L_{0},L_{1})\, |\,  
 \begin{array}{c}\lim_{s\ri-\infty}v(s,\cdot)=x\\ 
\lim_{s\ri+\infty}v(s,\cdot)=y\end{array} \right\}.$$

The following was  proved by by A. Floer in \cite{F1}, \cite{F3}  and by A. Floer, H. Hofer and D. Salamon
in \cite{FHS}:
\begin{theo}\label{floer}
(a) We have:
 $$\calm(L_{0}, L_{1})\, =\, \bigcup_{x,y\in L_{0}\cap 
L_{1}}\calm(x,y).$$
(b) For a
generic choice of  $J$ the spaces 
$\calm(x,y)$ are finite dimensional
manifolds with local dimension at $v\in\calm(x,y)$  given by the Maslov-Viterbo index  $\mu(v)$ (see
\cite{V} for the definition). 
\end{theo} 

 Denote by ${\call}^{0}(x,y)$ the zero-dimensional component of
${\call}(x,y)=\calm(x,y)/\real$. To define the   differential of $FC_{\bullet}$ we need
to prove that $\call^{0}(x,y)$ is finite. This is a consequence of Gromov's compactness for
holomorphic curves \cite{G} and was proved by A. Floer in \cite{F1} and Y.-G. Oh in \cite{Oh.}: 

\begin{theo}\label{gromov}  Suppose that $L$ is exact, or monotone with $N_{L}\geq 3$. 
Let $x,y\in L_{0}\cap L_{1}$ and $A> 0$ and let
$(v_{n})\subset\calm(x,y)$ be a sequence of solutions with constant index
$\mu(v_{n })=\mu_{0}\leq 2$. Then there exist a finite collection $(z_{i})_{i=0,\ldots,k}$
of  points in $L_{0}\cap L_{1}$ with $z_{0}=x$ and $z_{k}=y$,   some solutions 
$v^{i}\in\calm_{A}(z_{i-1},z_{i})$ for 
$i=1,\ldots k$ and some sequences of real
numbers $(\sigma_{n}^{i})_{n}$ for $i=1,\ldots,k $ such that for all $i=1,\ldots,k$
 a subsequence of
$v_{n}(s+\sigma_{n}^{i},t)$ converges  towards $v^{i}(s,t)$ 
in $\calc^{\infty}_{loc}$.

  Moreover, we have the relation
$$\sum_{i=1}^{k}\mu(v^{i})\, =\, \mu_{0}$$
\end{theo}

If $\mu_{0}=1$ in the statement above then necessarily $k=1$ due to the latter relation,
 so we immediately infer that the spaces $\call^{0}(x,y)$ are compact, which implies that
 they are
finite. This enables one to define the differential $\partial:FC_{\bullet}\ri FC_{\bullet}$ as:
$$\partial(x)\, = \, \sum_{y\in L_{0}\cap L_{1}}n(x,y)y,$$
where $n(x,y)=\#{\call}^{0}(x,y)$ $mod(2)$.
 In order to prove the relation $\partial\la{2}=0$, one has
to study the compactness of the $1$-dimensional component ${\call}^{1}(x,y)$ of ${\call}(x,y)$.
Using again Theorem \ref{gromov}  we find: 

\begin{theo}\label{compact} Denote by $\bar{\call}\la{1}(x,y)$ the union 
$${\cal L}^{1}(x,y)\, \cup\, \bigcup_{z\in L_{0}\cap
L_{1}}{\call}^{0}(x,z)\times{\call}^{0}(z,y),$$
endowed with the topology given by the convergence towards broken orbits which was
defined in Theorem \ref{gromov}. 

Then $\bar{\call}^{1}(x,y)$ is a compact 
$1$-dimensional manifold whose boundary is 
$\bigcup_{z\in L_{0}\cap L_{1}}{\call}^{0}(x,z)\times{\call}^{0}(z,y).$
\end{theo}

Note that the proof of the fact that $\bar{\call}\la{1}(x,y)$ is a manifold with boundary requires
a gluing argument as in \cite{F2}. Now the fact that a compact $1$-dimensional manifold has a
boundary of even cardinality immediately implies $\partial\la{2}=0$, proving thus that
$(FC_{\bullet}, \partial)$ is a complex. \\

 The following three subsections contain the proof of Theorem \ref{lift}. We start with: 

\subsection{Proof of Theorem \ref{lift} : Construction of the lifted Floer complex.}

A construction of a lifted Floer-type complex was already sketched in our previous work \cite{Da3}. 
The idea of constructing such a complex was suggested in \cite{BaCo}. 

 Consider the intersection points $L\cap \phi_{1}(L)$, 
viewed as points in $L$.
For two such points $x, y$, any holomorphic strip $v\in \calm(x,y)$ defines a path 
$\gamma:]-\infty,
+\infty[ \ri L$ which joins $x$ and $y$: 
$$\gamma(s) \, =\, v(s,0).$$
 We consider the obvious extension of $\gamma$ to $[-\infty,+\infty]$ keeping the same notation
for the extended path. Take the collection of intersection points and  the paths $\gamma$ as 
above, defined by the
strips $v$ which belong to the one-dimensional components of $\calm(x,y)$ 
(which correspond to the
zero dimensional components of $\call(x,y)$). 
 Denote by $\calc$ the collection of points and by $\Gamma$ the collection of
paths. 

Now start with the above collection $(\calc, \Gamma)$
 and fix a covering $p:\barl\ri L$.
 For any point $x\in \calc$ 
denote by $(x_{i})_{i\in I}$
the elements of the fiber $p\la{-1}(x)$. Consider all the lifts of the paths of $\Gamma$
 to the
covering space $\barl$. It is clear that for fixed points $x_{i}, y_{j}$ 
(where  $i,j\in
I$), the lifted space $\call\la{0}(x_{i},y_{j})$ is finite; its cardinality is obviously less then 
$\#\call\la{0}(p(x_{i}),p(y_{j}))$. Let $n(x_{i},y_{j})$ be 
the parity of this cardinality.
 On the free $\bf Z/2$-complex $C_{\bullet}\la{\barl}$ spanned by $\bigcup_{x\in
L\cap\phi_{1}(L)}p\la{-1}(x)$ one can  therefore define an
application $\partial\la{\barl}: C_{\bullet}\la{\barl}\ri C_{\bullet}\la{\barl}$ by the 
formula 
$$\partial\la{\barl}(x_{i})\, =\, \sum_{p(y_{j})=y\, \in\,  \calc}n(x_{i},y_{j})y_{j}.$$

The sum above is obviously well-defined since $\Gamma$ is finite and 
any path $\gamma\in \Gamma$ admits only one lifting starting from
$x_{i}$. We prove 

\begin{propo}\label{complex}
Under the hypothesis of Theorem \ref{lift} $(C_{\bullet}\la{\barl},\partial\la{\barl})$ is a complex.
\end{propo}

\noindent\underline{Proof}

   This is
equivalent to the fact that, given $x_{i}, y_{j}$, there is an even number of "broken paths"
joining them. Broken paths means liftings of  concatenations
$\gamma_{1}\ast\gamma_{2}$, where $\gamma_{i}\in \Gamma$ and $\gamma_{1}$ starts from $x$ and
ends into some point $z\in \calc$, while
$\gamma_{2}$ starts from the same point $z$ and  ends in $y$. 

Since the paths of $\Gamma$ define a complex (namely the Floer complex $FC_{\bullet}$), we know
that the number of broken paths joining $x$ and $y$ downstairs is even. Moreover, 
they represent
boundary points of a one-dimensional compact manifold, so they can naturally be dispatched 
in a
disjoint union of sets of two elements, corresponding to the boundary points of each 
component of the mentioned one-dimensional manifold. But in order to get the same property 
at the
level of the covering space $\barl$ one has to check that the two broken paths in such a pair admit
liftings to $\barl$ which have the same endpoints: 

\begin{propo}\label{pair} Let $\{\gamma_{1}\ast\gamma_{2}, \gamma_{1}'\ast\gamma_{2}'\}$ 
be a set of two broken
paths in $L$ as above. Then these broken paths are homotopic in $L$ with fixed endpoints. 
\end{propo}

\noindent\underline{Proof}

The concatenations  $\gamma_{1}\ast\gamma_{2}$ and $\gamma_{1}'\ast\gamma_{2}'$ have the same
starting point $x\in \calc$ and the same ending point $y\in \calc$. We use the following lemma which
was proved in \cite{Da3} (Lemma 3.16). 

\begin{lem}\label{homotopic} Let $(v_{n})_{n\in{\bf N}}$ be a sequence in $\calm(x,y)$, as in Theorem
\ref{gromov}.  Let $\gamma_{n} : [-\infty,+\infty]\ri L_{0}$ be  the path
defined by $\gamma_{n}(s)=v_{n}(s,0)$ extended with $x$ in $s=-\infty$ and with $y$ in
$s=+\infty$. For $i=1,\ldots, k$ let $\gamma^{i}:[-\infty,+\infty]\ri L_{0}$ be the
analogous paths defined by the holomorphic strips $v^{i}$. Then for $n$
large enough $\gamma_{n}$ and
$\gamma=\gamma^{1}\ast\gamma^{2}\ast\cdots\ast\gamma^{k}$ are homotopic in $L_{0}$. 
\end{lem}

Now the broken holomorphic strip $(v\la{1}, v\la{2})$ which defines $\gamma_{1}\ast\gamma_{2}$
corresponds to a boundary point in $\bar{\call}\la{1}(x,y)$, which means that it is the limit of a
sequence lying in a $2$-dimensional component of $\calm(x,y)$. The same is true for the broken path 
$\gamma_{1}'\ast\gamma_{2}'$.  Using the previous lemma we infer 
that there is some component of $\call\la{1}(x,y)$ such 
 $\gamma_{1}\ast\gamma_{2}$ and $\gamma_{1}'\ast\gamma_{2}'$ are respectively 
homotopic to  paths in $L_{0}$ 
defined by some elements in this component.  On the other hand all 
paths defined by the elements of
the same component of $\call\la{1}(x,y)$ are obviously homotopic, 
which finishes the proof of our
claim  \ref{pair}.

\hfill $\diamond$

The proof of Proposition \ref{complex} immediately follows, since the previous proposition implies that the
set of paths in $\Gamma$ whose lifts in $\barl$ join some fixed points $x_{i},\,  y_{j}$ is a
disjoint union of sets with two elements. 

\hfill $\diamond$

\vspace{.2in}

\noindent{\it Hamiltonian invariance} 

The usual way (see \cite{F2}, \cite{Oh.}) to prove that the homology 
of the Floer complex $FC_{\bullet}$ 
does not depend on the
(generic) choice of the Hamiltonian and of the almost complex structure, is to consider a
generic homotopy
$\Psi_{s}=(H_{s}, J_{s})$ between two fixed  couples $(H, J)$ 
and $(H', J')$ and to
use it to define a chain morphism between $FC_{\bullet}(H, J)
\ri FC_{\bullet}(H', J')$ which
induces an isomorphism in homology. More precisely, the homotopy $\Psi$
 is used to define moduli spaces
$\calm_{\Psi}(x,y)$ for $x\in\phi_{1}\la{H}(L)\cap L$ and 
$y\in\phi_{1}\la{H'}(L)\cap L$. The
 space $\calm_{\Psi}(x,y)$ is defined as follows:
 $$\left\{v:\real\times[0,1]\ri M\, 
\left|\begin{array}{c}\, 
\frac{\partial{v}}{\partial s}+J_{s}
\frac{\partial{v}}{\partial t}=0\\
\\ v(s,0)\in L_{0},\, \,
v(s,1)\in \phi_{1}\la{H_{s}}(L_{0}) \\ 
\\
\lim_{s\ri-\infty}v(s,t)=x, \,
\lim_{s\ri+\infty}v(s,t)=y\end{array} \, \right.\right\}.$$ 
For a generic choice of $\Psi$ these spaces are actually finite
dimensional manifolds, and the
morphism between the two complexes is defined by counting the number 
of elements (mod $2$) of their
$0$-dimensional components. The latter are proved to be finite by using a compactness result,
analogue to Theorem \ref{gromov}. The same result shows that the $1$-dimensional 
component of $\calm(x,y)$
can be completed to a compact boundary manifold such that the parity of 
the number of boundary
points is equivalent to the fact that the morphism defined by $\Psi$ 
commutes with the Floer
differentials.\\

The same argument as before can be used to get the invariance of the homology of 
$C_{\bullet}\la{\barl}$.
Given two collections of points and paths $(\calc,\Gamma)$ and $(\calc', \Gamma')$, as above, 
the
paths $v(s,0)\in \Gamma_{\Psi}$ defined by elements $v$ belonging to
$0$-dimensional components of $\calm_{\Psi}(x,y)$ define a
 morphism between the associated
lifted complexes $C_{\bullet}\la{\barl}$ and $(C_{\bullet}\la{\barl})'$.  
 Claiming that it is a chain morphism is
equivalent to  the fact that 
the number of broken paths in $\Gamma_{\Psi}$ admitting lifts which join fixed points
 $x_{i}$ and $y_{j}$ is even. As above, this is a consequence of Proposition \ref{pair}, adapted in
this new setting.

Finally, the arguments of \cite{F1} which show that the morphism defined by $\Psi$ between the two
usual Floer complexes induces an isomorphism at the homology level, can be used together with Proposition
\ref{pair} in the same way, in order
 to show that the homology of the lifted complex $C_{\bullet}\la{\barl}$ does not depend on
$(H_{t}, J)$ either. We will denote this homology by $FH\la{\barl}(L)$. \\

\subsection{Proof of Theorem \ref{lift} :  Computation of the lifted Floer homology  $FH\la{\barl}(L)$.} 

When the Lagrangian $L$ is weakly exact its Floer homology 
is isomorphic to the singular homology of
$L$ \cite{F2}. To prove this, one has to consider a Morse function $f:L\ri\real$ and a 
 particular Hamiltonian isotopy $\phi_{t}$
which maps $L$ into the
graph $L_{tf}\subset U(L)\subset M$, so that the Lagrangian intersections correspond to the
critical points of $f$.  The notations are those from \S1.3, in particular $U(L)$ is a
Weinstein tubular neighborhood of $L$. If $f$ is sufficiently $\calc\la{2}$-small
 one proves that 
 all the
holomorphic strips lie in $U(L)$; the contrary would imply - via Gromov compactness 
- the existence
of a nonconstant holomorphic disk with boundary in $L$ which is impossible for a weakly exact 
Lagrangian \cite{Oh}.
For a well chosen almost complex structure the canonical projection $U(L)\ri L$
maps the holomorphic strips onto the flow lines of a gradient vector field of $f$ 
with respect to a
generic Riemannian metric on $L$. Moreover, this projection defines a one-to-one 
correspondence between the
isolated  holomorphic strips and the gradient lines joining critical points of 
consecutive indices. This  
means that for these particular choices (which still satisfy
the genericity assumptions required to define Floer homology), the Floer complex is 
identical to a
Morse complex and the result follows.  

The latter property of the holomorphic strips above shows that in this particular case  the
collection $(\calc, \Gamma)$ is the one defined by the critical points of $f$ and the isolated
gradient lines which join them. So the lifted complex $C_{\bullet}\la{\barl}$ coincides 
with the
lifted Morse complex on $\barl$. The homology of the latter is the singular homology of
$\barl$ (recall that the stable manifolds of the gradient vector field associated to a 
Morse function
 yield a CW-decomposition of
$L$ whose lift to $\barl$ computes the homology of this covering space). Therefore we infer: 
$$FH\la{\barl}(L)\, \approx\, H_{\ast}(\barl,{\bf Z/2}).$$
 
In the monotone case, when one chooses the same particular Hamiltonian isotopy and 
almost complex
structure,  it is no longer true that the holomorphic strips lie in a Weinstein neighborhood
$U(L)$. But, as Y.-G. Oh pointed out in \cite{Oh}, the holomorphic strips which lie in 
$U(L)$ still
project onto the gradient lines and the isolated ones are in bijective correspondence 
to the gradient lines defining
the Morse complex. On the other hand, according to the same paper \cite{Oh} 
 an isolated holomorphic strip of finite energy which leaves $U(L)$  connects two critical 
points
$x,y$ of $f$ which satisfy 
$$Ind(x) -Ind(y) \, = \, 1-lN_{L},$$ for some positive integer $l$. 
  
In this particular case, the Floer complex can be graded by the Morse index. Therefore, 
given an integer $l$, the count  (mod $2$) of the isolated holomorphic strips satisfying 
the index relation
above defines for each integer $k$ a map $$\partial_{l} : FC_{k}\ri FC_{k-1+lN_{L}}.$$
Of course, $\partial_{l}$ vanishes for $l>\left[\frac{dim(L)+1}{N_{L}}\right]$.
The Floer differential $\partial : FC_{\bullet}\ri FC_{\bullet}$ writes
$$\partial\, =\, \partial_{0}+\partial_{1}+\cdots+\partial_{l}+\cdots.$$
Here $\partial_{0}$ is the Morse differential defined by the (projections on $L$ of the) 
 holomorphic strips which do not leave $U(L)$. 

From these data Y.-G. Oh \cite{Oh} and P. Biran \cite{Bi} inferred
 the existence of a spectral sequence which converges towards the Floer
homology and whose first page is built using the usual (Morse) homology of $L$. \\

Now, by definition, it is obvious that the differential $\partial\la{\barl}$ of the 
lifted complex 
 $C_{\bullet}\la{\barl}$ satisfies
 the same properties, namely:\\
1. It decomposes into a finite  sum
 $\partial\la{\barl}=\partial\la{\barl}_{0}+\partial\la{\barl}_{1}+\cdots$ \\
2.  For the grading given by the Morse index we have 
$\partial\la{\barl}_{l}:C_{k}\la{\barl}\ri C_{k-1+lN_{L}}\la{\barl} $\\
3.  The complex  $(C_{\bullet}\la{\barl},\partial\la{\barl}_{0})$ is identical 
to the lift to $\barl$ 
of a Morse complex on $L$. \\

These properties are sufficient for the existence of a spectral sequence analogous to the one 
defined by Biran and Oh. This  spectral sequence, denoted  $E_{r}\la{p,q}$,  is associated to an 
increasing filtration of a  complex
$\widetilde{C}_{\bullet}$ which is defined as follows: 

Denote by $A$ the subring ${\bf Z/2}[T\la{N_{L}}, T\la{-N_{L}}]$ 
of the Laurent polynomials with $\bf Z/2$-coefficients and by $A\la{kN_{L}}\subset A$ 
the subgroup 
${\bf Z/2}\cdot T\la{kN_{L}}$, for any integer $k$. Define the complex
 $$\widetilde{C}_{l}\, =\,
\bigoplus_{k\in{\bf Z}} C_{l-kN_{L}}\la{\barl}\otimes A\la{kN_{L}},$$
 endowed with the differential 
$$\tilde{d}\, =\, \partial_{0}\otimes Id\, + \partial_{1}\otimes T\la{-N_{L}}(\cdot)+ 
\partial_{2}\otimes T\la{-2N_{L}}(\cdot)+\cdots.$$
On this complex define a filtration ${\cal F}_{p}(\tilde{C}_{\bullet})$ by:

$$ {\cal F}_{p}(\widetilde{C}_{l})\, =\,
\bigoplus_{k\leq p} C_{l-kN_{L}}\la{\barl}\otimes A\la{kN_{L}},$$
It is easy to prove that the differential $\tilde{d}$ preserves the filtration and 
 that the homology of this complex is canonically  isomorphic to 
the Floer homology
$FH\la{\barl}(L)$ \cite{Bi}. Summarizing, the
 spectral sequence associated to these data has the following
properties:

\begin{theo}\label{spectral}  The  spectral sequence
 $\{E_{r}\la{p,q},d_{r}\}$ associated to the filtration $\cal F$
  converges to the lifted Floer homology  
$FH\la{\barl}(L)$ and  satisfies  the following properties 
(all the tensor products below are over $\bf Z/2$): \\
\begin{itemize}
\item
  $E_{0}\la{p,q}\, =\, C_{p+q-pN_{L}}\la{\barl}\otimes A\la{pN_{L}}$, $d_{0} =
[\partial\la{\barl}_{0}]\otimes Id$. 
\item $E_{1}\la{p,q}\, =\, H_{p+q-pN_{L}}(\barl, {\bf Z/2})\otimes A\la{pN_{L}}$, 
$d_{1}=[\partial\la{\barl}_{1}]\otimes T\la{-N_{L}}(\cdot)$, where 
$$[\partial_{1}\la{\barl}] : H_{p+q-pN_{L}}(\barl,{\bf Z/2}) \ri H_{p+q -1 -(p-1)N_{L}}
(\barl,{\bf Z/2})$$ is induced by $\partial_{1}\la{\barl}$. 
\item For every $r\geq 1$ $E_{r}\la{p,q}$ has the form $E_{r}\la{p,q}=V_{r}\la{p,q}
\otimes A\la{pN_{L}}$ with $d_{r}=\delta_{r}\otimes T\la{-rN_{L}}$, where $V_{r}\la{p,q}$ 
are vector spaces over $\bf Z/2$ and  $\delta_{r}:V_{r}
\la{p,q}\ri V_{r}\la{p-r,q+r-1}$ are homomorphisms
 defined for every $p,q$ and satisfying $\delta_{r}
\circ\delta_{r}=0$. 
Moreover: 
$$V_{r+1}\la{p,q}\, =\, \frac{Ker(\delta_{r}:V_{r}\la{p,q}\ri V_{r}\la{p-r,q+r-1})}
{Im(\delta_{r}:V_{r}\la{p+r,q-r+1}\ri V_{r}\la{p,q})},$$
and for $r=0,1$ we have $V_{0}\la{p,q}=C\la{\barl}_{p+q-pN_{L}}$, 
$V_{1}\la{p,q}=H_{p+q-pN_{L}}(\barl,{\bf Z/2})$, $\delta_{1}=[\partial\la{\barl}_{1}]$.
\item The spectral sequence collapses at page $\left[\frac{dim(L)+1}{N_{L}}\right]+1$ and 
for all $p\in {\bf Z}$, $\bigoplus_{q\in{\bf Z}}E_{\infty}\la{p,q}\, 
\approx\, FH\la{\barl}(L).$
\end{itemize}
\end{theo}

The proof of Theorem \ref{spectral} is analogous to the proof of Th. 5.2.A in \cite{Bi}. 
It is purely algebraic, and it applies to any  
graded complex whose differential (in our case
$\partial\la{\barl}$) satisfies  the  conditions 1-3 above.

 Remark that the above theorem immediately implies Theorem \ref{lift} for $L$ monotone. Indeed,  
 if $\delta_{1}=0$, then, according to Theorem \ref{spectral}, 
 $V_{2}\la{p,q}=V_{1}\la{p,q}=H_{p+q-pN_{L}}(\barl,{\bf Z/2})$, and 
$$\delta_{2}: H_{p+q-pN_{L}}(\barl,{\bf Z/2}) \ri H_{p+q-1-(p-1)N_{L}}(\barl,{\bf Z/2}).$$
Actually, the proof of \cite{Bi} shows that $\delta_{2}=[\partial_{2}\la{\barl}]$. 
Analogously, if
$\delta_{1}=\delta_{2}=\cdots=\delta_{l-1}=0$, then $\delta_{l}$ is defined on the homology of
$\barl$ (by $[\partial_{l}\la{\barl}]$) and its degree is $-1+lN_{L}$. Finally, if all the
$\delta_{i}$'s  vanish then the spectral sequence $E_{r}\la{p,q}$
 collapses at page $1$ and therefore, applying again
Theorem \ref{spectral}, we have
$$FH\la{\barl}(L)\, \approx\, H_{\ast}(\barl,{\bf Z/2}).$$
This finishes the proof of Theorem \ref{lift} for homologies with ${\bf Z/2}$-coefficients.  

\subsection{End of the proof of Theorem \ref{lift} : Change of coefficients}

A Lagrangian submanifold $L\subset M$ is called {\it relatively spin} 
 if there exists a class $st\in H\la{2}(M, {\bf Z/2})$ that restricts to the second
Stiefel-Whitney class $w_{2}(L)$ of $L$. For such submanifolds, also supposed to be
orientable, it was proved in \cite{FO3}
that the spaces of holomorphic strips $\calm(x,y)$ can be oriented and  that  under this
hypothesis
 the Floer complex $FC_{\bullet}$ can be defined over $\bf Z$-coefficients. 

In our case, when $L$ is orientable and relatively spin, 
the lifted complex $C_{\bullet}\la{\barl}$ is constructed using a 
collection of oriented paths $\Gamma$, which  enables us to define it over $\bf Z$. To
show that it is a $\bf Z$-complex whose homology 
only depends on $L\subset M$ and on the chosen
covering space $\barl$ one can use the same
proof as above. Analogously we get a spectral sequence $E_{r}\la{p,q}$ 
associated to a filtered complex as above.
But since short exact sequences over $\bf Z$ do not always split, we cannot infer as in Theorem 
\ref{spectral} that the lifted Floer homology is a direct sum of modules $E_{\infty}\la{p,q}$.
However, when the former vanishes, we know that $E_{\infty}\la{p,q}=0$ for all $p,q$. On the other
hand, the hypothesis of Theorem \ref{lift} implies that the spectral sequence collapses at page
$l+1$, and since $\delta_{m}=0$ for $m=1,\ldots, m-1$ (again by the hypothesis of Theorem
\ref{lift}), we find that $V_{l}\la{p,q}=H_{p+q-pN_{L}}(L,{\bf Z})$, whereas $V_{l}\la{p,q}=0$.
 Therefore,  the
complex $(H_{\ast}(\barl,{\bf Z}), \delta_{l})$ has to be acyclic.

 The proof of Theorem \ref{lift} is now complete.

\hfill $\diamond$

\begin{rema}\label{change}
If we chose the universal cover $\widetilde{L}$ as covering space, then
$C_{\bullet}\la{\widetilde{L}}$ can be seen as a free, finite dimensional complex over
${\bf Z/2}[\piu(L)]$ (resp. over ${\bf Z}[\piu(L)]$ when $L$ is orientable and 
 relatively spin). Of course
we can change the coefficients by tensoring it with any ${\bf Z/2}[\piu(L)]$-module $R$, 
for instance
with the Novikov ring associated to some morphism $u:\piu(L)\ri{\bf Z}$. 

In all these situations the homology of the lifted Floer
 complex $C_{\bullet}\la{\covl}$ is  related in the same
manner as for ${\bf Z/2}$-coefficients (resp. for integer coefficients) to the homology
of $\covl$ with coefficients in the new  ring. In particular, in the case of the Novikov ring,
 the latter is the Novikov homology $H_{\ast}(L,u)$ (for definition and related
properties, see for instance \cite{Da1}, \cite{Da2}). 
\end{rema}

\section{Applications}

In this section we prove our main results which 
we stated in \S1 and other applications of Theorem \ref{lift}.
\subsection{Aspherical Lagrangian submanifolds. Proof of 
 theorems \ref{main1} and \ref{coro1}}

The idea of the proofs is that, under the given hypothesis, 
the spectral sequence given by Theorem \ref{lift}
collapses at page $1$, which means that the lifted Floer homology associated to the
universal cover of $L$ 
is isomorphic to the singular
homology of $\widetilde{L}$. On the other hand, since $L$ is displaceable through 
a Hamiltonian isotopy, 
the lifted Floer 
homology vanishes. This is contradictory and therefore Theorem  \ref{lift} should 
not apply here. The only
possible reason for that  is
the fact that the Maslov number  $N_{L}$ is less than $3$.\\ 

\noindent\underline{ Proof of Theorem \ref{main1}}\\
(a) If $N_{L}\geq 3$ then we get the applications $\delta_{i}$ provided by 
Theorem \ref{lift}. But since $L$ is
aspherical $H_{i}(\widetilde{L})=0$ for $i\neq 0$, which implies that 
$\delta_{i}=0$ for all $i$,
therefore, according to Theorem \ref{lift}:
$$H(C_{\bullet}\la{\widetilde{L}})\, \approx\, H(\widetilde{L},{\bf Z/2}).$$
Since $L$ is displaceable, the left term vanishes, whereas
 the right term is not zero in degree $i=0$. This contradiction implies 
$N_{L}\leq 2$. 

We use the following well-known result \cite{Ar2} (see also \cite{Fu1}, Lemma 2.5):

\begin{propo}\label{orientable} If $L$ is orientable then   $N_{L}$ is even. The converse is true if 
$\piu(M)$ is trivial.
\end{propo}

 The conclusion of Theorem \ref{main1} follows. \\
\\ 
(b) Since $L$ is orientable, its Maslov number $N_{L}$ is even and 
therefore all the
applications $\delta_{i}$ provided by Theorem \ref{lift} have an odd degree. 
Suppose $N_{L}\geq
3$. The homology
of $\widetilde{L}$ is zero in odd degrees. This implies $\delta_{1}=0$ and
$E_{2}\la{p,q}=E_{1}\la{p,q}$ in the spectral sequence of Theorem \ref{spectral}. 
The same
argument shows that all the applications $\delta_{i}$ vanish,  
which implies that the
singular homology of $\widetilde{L}$ is zero. But this is impossible and therefore
$N_{L}=2$.\\
\\
(c) Denote by ${\cal J}_{reg}$ the generic set of compatible almost complex structures for
which the usual Floer complex $(FC_{\bullet}, \partial)$ is defined. Consider $J\in {\cal
J}_{reg}$ and denote 
$$\calm(M, J ; 2) \, =\, \{\, w:(D,\partial D)\ri (M,L)\, |\, \bar{\partial}_{J}w=0, \, \mu(w)=2\, \}. $$
 
By standard transversality results (\cite{DS}, see also \cite{BiCo2}, chap. 3) one gets
that for generic $J$,  $\calm(M, J; 2)$ is a manifold of dimension $n+2$. 
It is important to notice here
that a crucial point in the proof of the transverality is the fact that all
the disks in $\calm(M,J;2)$ are simple. This is a consequence of the monotonicity of $L$
and of a result of L. Lazzarini \cite{Laz} (see again \cite{BiCo2}). The monotonicity of
$L$ also implies that $\calm(M,J;2)$ is closed. Indeed all the holomorphic disks of this
manifold have the same area, so Gromov's compactness \cite{G} applies. On the other hand,
since $L$ is monotone and the disks have minimal Maslov number, no bubbling can occur. 

The unparametrized $J$-holomorphic disks of Maslov number $2$ passing through a given point
$p\in L$ can be identified with the preimage $ev\la{-1}(p)$ of an evaluation map 
$$\mbox{ev}\, :\, {\cal N}\, \ri\, L,$$ 
where ${\cal N}=(\calm(M,J;2)\times{\bfs}\la{1})/PSL(2,\real)$ and
 $PSL(2,\real)=Aut(D)$ acts on $(\calm(M,J;2)\times{\bfs}\la{1})$ by 
$$h\cdot(w,z)= (w\circ h, h\la{-1}(z)).$$
The evaluation map is given by $\mbox{ev}([w,z]) = w(z)$, for $[w,z]\in{\cal N}$. 

  The closed manifold ${\cal N}$ is $n$-dimensional and in particular for a generic $p$,
the preimage
$\mbox{ev}\la{-1}(p)$ is finite. Following the notations of \cite{Oh2}, denote by
$\Phi_{L}(p)$ the number of elements of $\mbox{ev}\la{-1}(p)$, modulo $2$. This number
does not depend on the choice of the regular value $p$: it is the mod-$2$ degree of the
evaluation map. Y-G. Oh shows in \cite{Oh2} that in order to define the Floer homology
$FH(L_{0},L_{1})$ in the case where $N_{L_{i}}=2$ one needs the hypothesis
$\Phi_{L_{0}}+\Phi_{L_{1}}=0$. This comes from the fact that the $2$-dimensional component
of the trajectory spaces $\calm(x,y)$ can be compactified by adding the broken trajectories
{\it and} the holomorphic disks of Maslov index $2$ with boundary in one of the $L_{i}$'s,
 passing through the intersection
points; the latter occur as bubbles of sequencences in $\calm(x,y)$ (when $x=y$). 
If their number is
even then the relation $\partial\la{2}=0$ is still valid. 

Note that Oh also shows that
when $L_{1}$ is a Hamiltonian deformation of $L_{0}$ then the above relation 
is satisfied and
therefore the Floer homology $FH(L)$ can be defined.  

But for the definition of the lifted Floer
homology $FH\la{\covl}(L)$ this  is no longer sufficient as it can be seen in 
 the example of two circles in
$\bf C$ intersecting in two points. Recall that the lifted Floer complex was defined (for
$N_{L}\geq 3$) using the
paths $w(s,0)$ defined by the isolated holomorphic strips of $\calm(x,y)$.
To define it for $N_{L}=2$ (with $\bf Z/2$-coefficients) one 
needs for any $x\in
L\cap L_{1}, \, \, L_{1}=\phi_{1}(L)$ and for any  homotopy class $g\in\piu(L)$ 
an even number of broken isolated trajectories $w$ from $x$ to $x$ 
whose associated paths $w(s,0)\subset L$ define a loop in the class $g$. This is clearly not
true in the case of the two circles, as the two broken paths lie in different homotopy
classes.

 In the general situation, we claim that the Floer homology $FH\la{\covl}(L)$
can be defined when for
any $g \in \piu(L)$ the number of $J$-holomorphic
 disks passing through a generic $p\in L$ and whose
boundary realize $g$ is even. The same arguments as above 
(fixing the homotopy class of the boundary in
the definitions of $\calm$, ${\cal N}$) show that the parity of 
this number does not depend on the generic choice of $p$.
  Denote it  by $\Phi_{g, L} \in\{0,1\}$; clearly 
$$(1)\gol \Phi_{L}\, =\, \sum_{g\in\piu(L)}\Phi_{g, L}\gol \mbox{mod}\, 2 .$$

Let us prove our claim. As explained above, 
it can happen that in order to compactify an 
one-dimensional connected component of $\call\la{1}(x,x)$ (where $x\in
L\cap\phi_{1}(L)$), 
one has
to add a broken trajectory and a holomorphic disk $D$ with boundary in 
$L$ or $\phi_{1}(L)$. 
The homotopy class of the boundary $\partial D\subset L$ is determined by
 the  paths $w(s,0)$ defined by the holomorphic
strips corresponding to the elements in $\call\la{1}(x,x)$. When the boundary contains a
holomorphic disk $(D,\partial D)\subset (M,L_{1})$, then the loops defined in $L$ by the
trajectories of $\call(x,x)$ are necessarily contractible in $L$. Therefore, counting modulo
$2$, we have $\Phi_{g,L}$ broken trajectories in the class $g$ when $g\neq 0$ and
$\Phi_{0,L}+\Phi_{L_{1}}$ for $g=0$. Using the fact that $\Phi_{L_{1}}=\Phi_{L}$
and the relation $(1)$ we
infer that the lifted Floer complex is defined provided that $\Phi_{g,L}=0$ for any
$g\in\piu(L)$, $g\neq 0$. 

Now we are able to finish our proof. If for any non zero $g\in\piu(L)$ there is 
some $p\in L$ such that there is no holomorphic disk with
boundary in the class $g$, passing through $p$, then $\Phi_{g,L}=0$ for any $g\neq 0$,
 so $FH\la{\covl}(L)$ is defined.
But this leads to a contradiction like in the proof of Theorem \ref{main1}.(a). So, the proof is finished
for $J\in {\cal J}_{reg}$. For an arbitrary $J$, take a sequence $J_{n}\in {\cal J}_{reg}$ which
converges towards $J$. Fix $p\in L$ and consider $J_{n}$-holomorphic disks $w_{n}:(D,\partial
D)\ri (M,L)$ such that $\mu(w_{n}) =2$ and $p\in w_{n}(\partial D)$. Using again Gromov's
compactness \cite{G} we find that $w_{n}$ converges towards a $J$-holomorphic disk whose boundary passes
through $p$. There is no bubbling here because of the monotonicity of $L$ and of the fact that the
Maslov index is minimal. The boundary of the limit is not trivial in $\piu(L)$ by an argument which is
similar to the one in the proof of Lemma \ref{homotopic} (see \cite{Da3}, Lemma 3.16). 

The proof of Theorem \ref{main1} is now complete.

\hfill $\diamond$

\vspace{.2in}

 Using a quite similar argument one can prove the following version of
Theorem \ref{main1}:
\begin{theo}\label{main1gen} Let $M$ be a monotone symplectic manifold which has the
property that any compact subset is displaceable through a Hamiltonian isotopy and let
$L\subset M$ be a monotone Lagrangian submanifold. \\
(a) If for some integer $k \geq 1$ we have 
$H_{i}(\widetilde{L}, {\bf Z/2})=0$ for $i > k$ then $$N_{L}\in [1, k+1]$$ 

(b) Suppose that $L$ is orientable. If $H_{\ast}(\widetilde{L},{\bf Z/2})$ is of finite 
dimension over $\bf Z/2$ and the
Euler characteristic
$$\chi =\sum_{i=0}\la{n}(-1)\la{i}dim(H_{i}(\widetilde{L},{\bf Z/2}))$$ does not vanish,
then $N_{L}=2$. Moreover, for any 
almost complex structure $J$ which is compatible with the symplectic form we have that 
 through every $p\in L$ passes at least 
a $J$-holomorphic disk $w:(D,\partial D)\ri (M,L)$  
whose Maslov
index equals $2$ and whose boundary is not trivial in $\piu(L)$. \\

\end{theo}

\noindent\underline{Proof}

(a) If $N_{L}\geq k+2\geq 3$ the lifted Floer homology $FH\la{\widetilde{L}}$ 
is well defined.
Since the degree of the applications $\delta_{i}$ is greater of equal to 
$-1+N_{L}\geq
k+1$, all these application vanish and therefore $FH\la{\widetilde{L}}$ 
is isomorphic to
the singular homology of $\widetilde{L}$. On the other hand 
$FH\la{\widetilde{L}}=0$, as
$L$ is displaceable, which is absurd.

(b) To prove this statement one has to look at the proof of Th. 5.2.A in 
\cite{Bi}. It is shown that the
vector spaces $V_{r}\la{p,q}$ in the statement of Theorem \ref{spectral} satisfy
$V_{r}\la{p+1,q}=V_{r}\la{p,q+1-N_{L}}$, and the applications 
$\delta_{r}\la{p,q}$ have the same
property. Therefore, for $p, q$ fixed, $\delta_{r}$ is a differential 
on the complex
$(V_{r}\la{p,q+k(1-rN_{L})})_{k\in{\bf Z}}$. This complex is finite 
(since this assertion  is true for $r=1$) and
its homology is $(V_{r+1}\la{p,q+k(1-rN_{L})})_{k}$, according to Theorem \ref{spectral}.
 Now
fix $p\in {\bf Z}$ and consider the Euler characteristic:
$$\chi_{r}\, =\, \sum_{q\in{\bf Z}}(-1)\la{q}dim(V_{r}\la{p,q}).$$
We show that $\chi_{r}$ does not depend on $r$. Fix a negative odd number $m$.
 It is quite clear that 
$$\chi_{r}=\sum_{l=m+1}\la{0}(-1)\la{l}\sum_{k\in{\bf Z}}(-1)\la{k}
dim(V_{r}\la{p,l+km}).$$
Indeed we have just changed the order of the summands in the
 writing of $\chi_{r}$. Applying this for $m=1-rN_{L}$ 
(which is odd, since $L$ is orientable, by Proposition \ref{orientable}), we get:
$$\chi_{r}\, =\, \sum_{q=2-rN_{L}}\la{0}(-1)\la{q}\chi_{r}\la{q},$$
where $\chi_{r}\la{q}$ is the Euler characteristic of the complex
 $(V_{r}\la{p,q+k(1-rN_{L})})_{k}$. The Euler characteristic of the homology 
is the same, so one can write 
$$\chi_{r}\la{q}\, =\, \sum_{k\in{\bf Z}}(-1)\la{k}
dim(V_{r+1}\la{p,q+k(1-rN_{L})}),$$
therefore we get
$$\chi_{r} =\sum_{q=2-rN_{L}}\la{0}(-1)\la{q}\sum_{k\in{\bf Z}}(-1)
\la{k}dim(V_{r+1}\la{p,q+k(1-rN_{L})})=\chi_{r+1}.$$
The latter equality is obtained by applying the property above for 
$m=1-rN_{L}$ and $r+1$ instead of $r$. 

So $\chi_{r}$ is independent of $r$. On the other hand, 
according to Theorem \ref{spectral}, we have $$\chi_{1}=\chi(H_{\ast}(\widetilde{L},
{\bf Z/2}))\, \neq\, 0$$
and since the spectral sequence collapses and its limit 
is zero we also have $\chi_{r}=0$ for $r$ sufficiently large. Therefore the
 lifted Floer complex can not be defined, which means that $N_{L}=2$.

The proof of the existence of a $J$-holomorphic disk of Maslov index $2$ passing through a
given $p\in L$ is similar to the proof of Theorem\ref{main1}.(a). 

\hfill $\diamond$

\vspace{.2in}

\noindent\underline{Proof of Theorem \ref{coro1}}

We use the following result of P. Biran (\cite{Bi}, Prop. 4.1.A):
\begin{propo}\label{biran} If $L\subset {\bf CP}\la{n}\times W$ 
is Lagrangian monotone then there is a circle bundle 
$\Gamma_{L}\ri L$ such that $\Gamma_{L}$ admits a monotone Lagrangian 
embedding into ${\bf C}\la{n+1}\times W$ and moreover $N_{\Gamma_{L}}=N_{L}$. 
\end{propo}

If $L$ is aspherical then it is easy to see that $\Gamma_{L}$ 
is aspherical, too. Since it is displaceable through a Hamiltonian isotopy, 
the result follows by Theorem \ref{main1}. 

\hfill $\diamond$

\subsection{Lagrangian submanifolds with maximal Maslov number. 
Proof  of  theorems \ref{main2} and \ref{projective}}

\noindent\underline{Proof of Theorem \ref{main2}}

(a) As $N_{L}=n+1\geq 3$ the lifted Floer complex is well defined. 
 We know that its homolgy $FH\la{\barl}(L)$ vanishes because $L$ is
displaceable. 
On the other hand, we know by Theorem \ref{spectral} that the spectral sequence
$\{E_{r}\la{p,q}\}$ which converges towards this homology, collapses at page $2$.
So,  according to Theorem \ref{spectral} we have 
$$0\, =\, E_{2}\la{p,q}\, =\, \frac{Ker\left([\partial_{1}\la{\barl}]:E_{1}\la{p,q}\ri
E_{1}\la{p-1,q}\right)}{Im\left([\partial_{1}\la{\barl}]:E_{1}\la{p+1,q}\ri
E_{1}\la{p,q}\right)},$$
which for $q=p-pN_{L}+i$ gives 
$$0=E_{2}\la{p,p-pN_{L}+i}\, =\, \frac{Ker\left([\partial_{1}\la{\barl}]:
H_{i}(\barl,{\bf Z/2})\ri H_{n+i}(\barl,{\bf Z/2})\right)}
{Im\left([\partial_{1}\la{\barl}]:
H_{i-n}(\barl,{\bf Z/2})\ri H_{i}(\barl,{\bf Z/2})\right)}$$

 Applying this equality for  $i=1, \ldots n-1$ we find that $\barl$ is a
 ${\bf Z/2}$-homology sphere for any covering space $\barl$. In particular, for $\barl=L$ we have  
$H\la{i}(L,{\bf Z/2})=0$, for $i=1,2$
 and therefore we infer 
that $L$ is  spin. So $dim(L)=N_{L}-1$ is odd (by Proposition \ref{orientable}, since $L$ is oriented)
and the whole theory works for $\bf Z$-coefficients. The same argument shows then that 
 any covering $\barl$ is a
 $\bf Z$-homology sphere. Let us prove that $L$ is also simply connected. If not, take 
a non zero element $g\in \piu(L)$ and consider the Abelian subgroup $G=<g>$ and 
the associated covering space $\barl$. Therefore $H_{1}(\barl, {\bf Z})=G$ which contradicts the fact
that $\barl$ is a $\bf Z$-homology sphere.  Finally,  using the (proofs of the) Poincar\'e
conjecture, (S. Smale, M. Freedman, G. Perelman) we infer that  $L$ is homeomorphic to ${\bfs}\la{n}$,
as claimed. \\

(b) As above, the lifted Floer complex is defined, its homology vanishes, and the
spectral sequence converging to it collapses at page $2$. Consider an arbitrary
 covering
$\barl\ri L$. We get, as in the proof of a) for $i=1, \ldots n$: 
$$E_{2}\la{p,p-pN_{L}+i}\, =\, \frac{Ker\left([\partial_{1}\la{\barl}]:
H_{i}(\barl,{\bf Z/2})\ri H_{i+n-1}(\barl,{\bf Z/2})\right)}
{Im\left([\partial_{1}\la{\barl}]:
H_{-n+i+1}(\barl,{\bf Z/2})\ri H_{i}(\barl,{\bf Z/2})\right)}.$$
We infer that $$ (1) \gol H_{i}(\barl,{\bf Z/2}) = 0$$ for $i=2, \ldots n-2$ and
$$(2) \gol [\partial_{1}\la{\barl}]:  H_{1}(\barl,{\bf Z/2}) \, \approx\,
H_{n}(\barl, {\bf Z/2}).$$
$$(3) \gol [\partial_{1}\la{\barl}]:  H_{0}(\barl,{\bf Z/2}) \, \approx\,
H_{n-1}(\barl, {\bf Z/2}).$$

Suppose first that $L$ is not orientable. Choosing $\barl=\widetilde{L}$ we 
get $H_{n}(\widetilde{L})=0$, so 
$\widetilde{L}$ is
not compact, and therefore $\piu(L)$ is infinite.  We know that 
$H_{1}(L, {\bf Z/2})\approx
H_{n}(L,{\bf Z/2})={\bf Z/2Z} $. Take an 
element $g\in \piu(L)$ which is not
in the kernel of the Hurewicz morphism $\piu(L)\ri H_{1}(L,{\bf Z/2})$ 
and consider the covering
$\barl\ri L$ associated to the Abelian subgroup $G=<g>$. If $g$ is of
 finite order - which has to be even -, the covering
$\barl\ri L$ is infinite (since $\piu(L)$ is infinite), 
therefore $H_{n}(\barl, {\bf Z/2})=0$ and 
 $H_{1}(\barl, {\bf Z/2})=0$, by $(2)$. On the other hand 
$H_{1}(\barl, {\bf Z/2})=G\otimes_{\bf
Z}{\bf Z/2}\neq 0$, since $G$ is a cyclic group ${\bf Z/2l}$ due to the fact that  $g$ has en even order.
 This contradiction implies 
 that $G$ is an infinite cyclic
group. In particular $H_{1}(\barl, {\bf Z/2})\neq 0$, and by $(2)$ 
$H_{n}(\barl, {\bf Z/2})\neq 0$,
which means that $\barl$ is compact and $G$ has finite index in $\piu(L)$, as claimed. 

If $M$ is an exact symplectic manifold then Theorem \ref{Gromov} applies and 
$H\la{1}(L,\real)\neq 0$, a
non-zero class being given by the restriction to $L$ 
of a primitive of the symplectic form. Therefore
the first Betti number of $L$ is not zero. Consider a non-vanishing morphism 
$u:\piu(L)\ri\bf Z$ and denote by $K$ its kernel. We show that $K$ is finite, of odd order. 

Indeed, take an element $t\in \piu(L)$ such that $u(t)=1$ and consider the covering
$\barl\ri L$ corresponding to the infinite cyclic subgroup $G$ spanned by $t$. As above, using
the relation $(2)$ we infer that $\barl\ri L$ is a finite covering. It is easy to see that
different elements of $K=Ker(u)$ lie in different classes of the quotient $\piu(L)/G$ and therefore
$K$ is finite. Moreover, if we suppose that an element $g\in K$ has even order, then, taking
as $\barl$ the covering associated to $G=<g>$, we find as above that $H_{1}(\barl,{\bf
Z/2})\neq 0$, whereas $H_{n}(\barl,{\bf
Z/2})= 0$ since the covering is infinite. This contradicts the relation $(2)$. \\

Let us consider now the case where  $L$ is orientable (so $n$ is even). Since
it is also spin (by $(1)$), the relations $(1)$ and $(2)$ are valid for integer
coefficients. Taking $\barl = L$, we get $H_{1}(L,{\bf Z})=\bf Z$, which implies
$H\la{1}(L,{\bf Z})={\bf Z}$. Therefore, as above there is an exact sequence of groups
$$0\ri K\ri\piu(L)\ri {\bf Z}\ri 0,$$ where $K$ is finite. Now, using the same argument as
above for the covering  $\barl$ associated to $<g>$, for an arbitrary $g\in K$, we find that
 $H_{1}(\barl,{\bf
Z})\neq 0$ unless $g$ is the identity. Using the relation $(2)$, this implies that
$K=\{1\}$, so $\piu(L)\approx {\bf Z}$, as claimed. 

For 
$dim(L)\geq 6$ F. Latour and A. Pajitnov independently established an algebraic criterion for
the existence of a fibration of $L$ over the circle \cite{Lat}, \cite{Paj}. For
$\piu(L)={\bf Z}$ we get (see for instance \cite{Da2}): 

\begin{theo}\label{fibration} When $n=dim(L)\geq 6$ and $\piu(L)=\bf Z$, then there exists a
fibration $f:L\ri\bfs\la{1}$ if and only if the Novikov homology $H_{\ast}(L;u)$ vanishes,
where $u=[f\la{\ast}(d\theta)]\in H\la{1}(L,{\bf Z})$. 
\end{theo}

According to Theorem \ref{lift} (see Remark \ref{change})
 the relations $(1)$ and $(2)$ are also valid for the Novikov
homology with respect to any $1$-cohomology class $u$. On the other hand, for any $u\neq 0$
one can show that $H_{0}(L;u)=0$ and $H_{n}(L;u)=0$ (see for instance \cite{Da1}). Using
$(1)$, $(2)$ and $(3)$ we find that for $u=id_{\bf Z}$, $H_{\ast}(L,u)=0$ and therefore,
$L$ admits a
fibration over the circle, by Theorem \ref{fibration}. Denote by $F$ a fiber of this fibration. We
know that $\piu(F)=Ker(u)$, so $F$ is simply connected. We also have that $\widetilde{L}$ is
diffeomorphic to $F\times \real$ and in particular $$H_{\ast}(\covl,{\bf Z})\approx
H_{\ast}(F,{\bf Z}).$$ We infer from $(2)$ that $F$ is a simply connected homology sphere,
therefore it is homeomorphic to the standard $(n-1)$-sphere, using Poincar\'e.

The proof of Theorem \ref{main2} is now complete.

\hfill $\diamond$

\vspace{.2in} 

\noindent\underline{Proof of Theorem \ref{projective}}

(a) The vanishing of the first homology group of $L$ implies that $L$ is monotone with $N_{L}= 2(n+1)$.
Applying Proposition \ref{biran} we get a monotone Lagrangian submanifold $\Gamma_{L}\subset {\bf C}\la{n+1}\times X$
which has the same Maslov number $N_{\Gamma_{L}}=N_{L}=dim(\Gamma_{L})+1$.  Moreover $\Gamma_{L}\ri L$ is a
circle fibration. 

The submanifold $\Gamma_{L}$ satisfies  the hypothesis of Th. \ref{main2}.(b) and therefore
$\piu(\Gamma_{L})\approx {\bf Z}$. Moreover, considering the lifted Floer complex associated to the
universal cover $\widetilde{\Gamma}_{L}$, the relations $(1)$, $(2)$ and $(3)$ from the proof of 
Theorem~\ref{main2} imply that $H_{i}(\widetilde{\Gamma}_{L})=0$ for $i\neq 0, 2n+1$, and 
$H_{2n+1}(\widetilde{\Gamma}_{L})\approx {\bf Z}$. Therefore, using Hurewicz's isomorphism
 the first $2n+1$ homotopy groups of
$\Gamma_{L}$ are isomorphic to the corresponding homotopy groups of $\bfs\la{2n+1}$. 
  From the long exact sequence of the fibration $\Gamma_{L}\ri L$ we infer that 
$$\pi_{i}(L)\, \approx\, \pi_{i}(\bfs\la{2n+1}), \, \, \mbox{for}\, \, \, i=2,\ldots 2n+1,$$
and, since $\pi_{1}(\Gamma_{L})$ is $\bf Z$,  $\piu(L)$ is Abelian. But as $H_{1}(L,{\bf Z})=0$, $L$ is
simply connected. Therefore $L$ is homeomorphic to $\bfs\la{2n+1}$ according to (the proof of)
 the Poincar\'e
conjecture. \\
\\
(b) As above $L$ is monotone with $N_{L}=2(n+1)$. Again, we can apply Proposition \ref{biran} and we get a
monotone Lagrangian in $\bfc\la{n+1}\times {\bf CP}\la{n}$  which is a circle fibration over $L$ and has
the same Maslov number. We can therefore use Th.\ref{main2}.(a) and infer that $\Gamma_{L}$ is homeomorphic
to the $(2n+1)$-sphere, as claimed. In particular $L$ is simply connected.\\
\\
(c)  As pointed out in \cite{Bi}, one can easily see that $L$ is monotone with $N_{L}=n+1$.
 We consider as
above the Lagrangian submanifold $\Gamma_{L}\subset \bfc\la{n+1}$ which has the same Maslov number.
 We can
use Theorem \ref{lift} with integer coefficients because $\Gamma_{L}$ is spin 
(using again Theorem \ref{lift} with 
mod-$2$
coefficients) and orientable (by Proposition \ref{orientable}, since $N_{\Gamma_{L}}$ is even). As in the proof 
of (a) above we get that
$\piu(\Gamma_{L})={\bf Z}$ and  the
universal cover $\widetilde{\Gamma}_{L}$ has the same homotopy groups $\pi_{i}$ as $\bfs\la{n}$, for
$i=1,\ldots n$. Using the long exact sequence of the fibration $\Gamma_{L}\ri L$ 
we find $\pi_{n}(L)={\bf
Z}$, $\pi_{i}(L)=0$ for $i=3,\ldots n-1$ and we have an exact sequence:
$$0\, \ri\, \pi_{2}(L)\, \ri\, \pi_{1}(\bfs\la{1})={\bf Z}\, \ri\, \piu(\Gamma_{L})={\bf Z}\, \ri\,
\piu(L)\, \ri\, 0.$$
It follows that $\piu(L)$ is cyclic Abelian and, since $H_{1}(L,{\bf Z})$ is $2$-torsioned, 
it follows that
$\piu(L)={\bf Z/2}$. From the exact sequence we infer then that $\pi_{2}(L)=0$. So the universal cover of
$L$ is a homotopy sphere and therefore it is homeomorphic to the  $n$-sphere. 

\hfill $\diamond$

 \subsection{Lagrangian submanifolds in the cotangent bundle. Proof of Theorem \ref{sphere}}

Let $K\subset M$ be a Lagrangian submanifold. If $L\subset T\la{\ast}K$ is Lagrangian then by 
Darboux' Theorem it
follows that $L$ admits a Lagrangian embedding into $M$. We need the following result:

\begin{propo}\label{monotone} Suppose that $K\subset M$ is monotone and that
 $L\subset T\la{\ast}K$ is exact with vanishing Maslov class. Then $L\subset
M$ is also monotone. Moreover, if the morphism $\piu(L)\ri\piu(K)$ induced by the projection is surjective, then
$N_{L}=N_{K}$. 
\end{propo}

\noindent\underline{Proof}  

Since $L$ is exact it is easy to see that (using the notations of the first section) we have 
 $$I_{\omega}\la{L\hookrightarrow M }=I_{\omega}\la{K\hookrightarrow
M}\circ p,$$ where $p:\pi_{2}(M,L)\ri\pi_{2}(M,T\la{\ast}K)\approx\pi_{2}(M,K)$ is the canonical morphism. 

 It is also known and not very hard to prove (see  \cite{Hol}, Chap.I, Prop. A.3.3) that when $L$ has a 
vanishing Maslov class
then (again with the  notations of \S1) we have 
$$I_{\mu}\la{L\hookrightarrow M }=I_{\mu}\la{K\hookrightarrow
M}\circ p.$$ 
Therefore $L\subset M$ is monotone. If $\piu(L)\ri\piu(K)$ is surjective,
 then $p$ is also an epimorphism and the
conclusion follows.

\hfill $\diamond$.  

\vspace{.2in}

\noindent \underline{Proof of Theorem \ref{sphere}.(a)} Let $L\subset T\la{\ast}\bfs\la{2k+1}$ be exact
Lagrangian  with vanishing
Maslov class.  As pointed out in \S1.2, 
the application 
$$z\, \mapsto\, ([z], \bar{z})$$
defines a monotone Lagrangian embedding of $\bfs\la{2k+1}$ into ${\bf CP}\la{k}
\times{\bf C}\la{k+1}$. It follows by
Proposition \ref{monotone} that $L$ admits a monotone Lagrangian embedding into ${\bf CP}\la{k}
\times{\bf C}\la{k+1}$ of Maslov
number $N_{L}=N_{\bfs\la{2k+1}}=2k+2$. The statement (a) of  Theorem \ref{main2} 
implies the desired result. \\

\noindent \underline{Proof of Theorem \ref{sphere}.(b)} Let 
$L\subset T\la{\ast}K$ be an exact Lagrangian  submanifold with vanishing
Maslov class. There is a finite cover $\barl$ of $L$ which admits an 
exact Lagrangian embedding into $T\la{\ast}\widetilde{K}=T\la{\ast}\bfs\la{2k+1}$ (see \cite{Da3}, Lemma
3.5). Moreover, we have a commutative diagram: 
$$\begin{array}{ccc}
\bar{L}&\, \, \ri\, \, &T\la{\ast}\bfs\la{2k+1}\\
\\
\downarrow\, &\, &\downarrow\\
\\
L&\, \, \ri\, \, &T\la{\ast}K
\end{array} $$ 
From this diagram one immediately infers that the Maslov class of $\barl$ vanishes. Therefore, we can
apply the point (a) which asserts that $\barl$ is homeomorphic to $\bfs\la{2k+1}$. 

When $K ={\bf RP}\la{2k+1}$ we show that 
 $\piu(L)\ri\piu({\bf RP}\la{2k+1})$ is  an epimorphism. If not,  we can lift
$L$ to an exact Lagrangian embedding (of $L$) into $T\la{\ast}\bfs\la{2k+1}$, 
which still has vanishing Maslov
 class. 
Using the statement (a) we find that $L$ is homeomorphic to a sphere. Since ${\bf RP}\la{2k+1}$ admits a
Lagrangian embedding into ${\bf CP}\la{2k+1}$, the same is true for $L$. 
 But this contradicts
 a well-known result of P. Biran and K. Cieliebak asserting that 
there is no Lagrangian sphere in ${\bf CP}\la{n}$ (\cite{BiC}, Theorem A).
 So $\piu(L)\ri\piu({\bf RP}\la{2k+1})$ is surjective. In this case, using Lemma 3.5 
of \cite{Da3} as above we get that $\barl=\covl$ is a double covering, homeomorphic to $\bfs\la{2k+1}$. 

\hfill $\diamond$

\vspace{.2in}

\noindent {\bf Acknowledgements} I thank Alexandru Oancea for our valuable discussions on the subject.
 I am grateful to Alexandru Oancea and Fran\c{c}ois Laudenbach, who
carefully read the manuscript and made several corrections and useful suggestions.

\end{document}